\documentclass[12pt]{amsart}
\usepackage{amssymb}
\usepackage{amsmath} 
\usepackage{amsthm} 
\usepackage{calc} 
\usepackage{verbatim} 
\usepackage[dvips]{graphics}
\newtheorem{thm}{Theorem}[section]
\newtheorem{lem}[thm]{Lemma}

\newtheorem{cor}[thm]{Corollary}

\newtheorem{prop}[thm]{Proposition}

\newcounter{probno}
\stepcounter{probno}

\newcommand{\proofof}[1]{\noindent{\it Proof of #1. }}
\renewcommand{\qed}{\hfill$\Box$ \par \medskip}

\newcommand{\ignore}[1]{}

\newcommand{\Area}{\mathrm{Area}\,}

\newcommand{\self}{\circlearrowleft}

\newcommand{\eqdef}{:=}

\newcommand{\id}{\mathrm{id}}

\newcommand{\R}{\mathbf{R}}

\newcommand{\C}{\mathbf{C}}
\newcommand{\cp}{\mathbf{P}}
\newcommand{\Z}{\mathbf{Z}}
\newcommand{\N}{\mathbf{N}}

\newcommand{\Q}{\mathbf{Q}}

\newcommand{\supp}{\mathrm{supp}\,}

\renewcommand{\eqdef}{:=}

\newcommand{\basin}{\mathcal{B}}

\title{Real Dynamics of a Family of Plane Birational Maps:  Trapping 
       Regions and Entropy Zero}

\date{\today}

\author{Eric Bedford \& Jeffrey Diller}
\address{Department of Mathematics\\
         Indiana University\\
         Bloomington, IN, IN 47405}
\email{bedford@indiana.edu}
\address{Department of Mathematics\\
         University of Notre Dame\\
         Notre Dame, IN 46556}
\email{diller.1@nd.edu}

\thanks{The first author is partially supported by the National Science Foundation}
\begin{document}

\maketitle
\markboth{\today}{\today}

\section{Introduction}
\label{intro}
We consider dynamics of the one parameter family of birational maps
\begin{equation}
\label{opf}
f = f_a:(x,y)\mapsto\left( y{x+a\over x-1},x+a-1\right).
\end{equation}
This family was introduced and studied by Abarenkova, Angl\`es d'Auriac, 
Boukraa, Hassani, and Maillard, with results published in [A1--7].  We 
consider here real (as opposed to complex) dynamics, treating $f_a$ as
a self-map of $\R^2$, and we restrict our attention to parameters $a>1$.  In 
order to discuss our main results, we let $\basin^+,\basin^-\subset \R^2$ be 
the sets of points with orbits diverging locally uniformly to infinity in 
forward/backward time, and we take $K\subset\R^2$ to be the set of points 
$p\in\R^2$ whose full orbits $(f^n(p))_{n\in\Z}$ are bounded.  

\begin{figure}
\begin{center}
\resizebox{3in}{3in}{\includegraphics{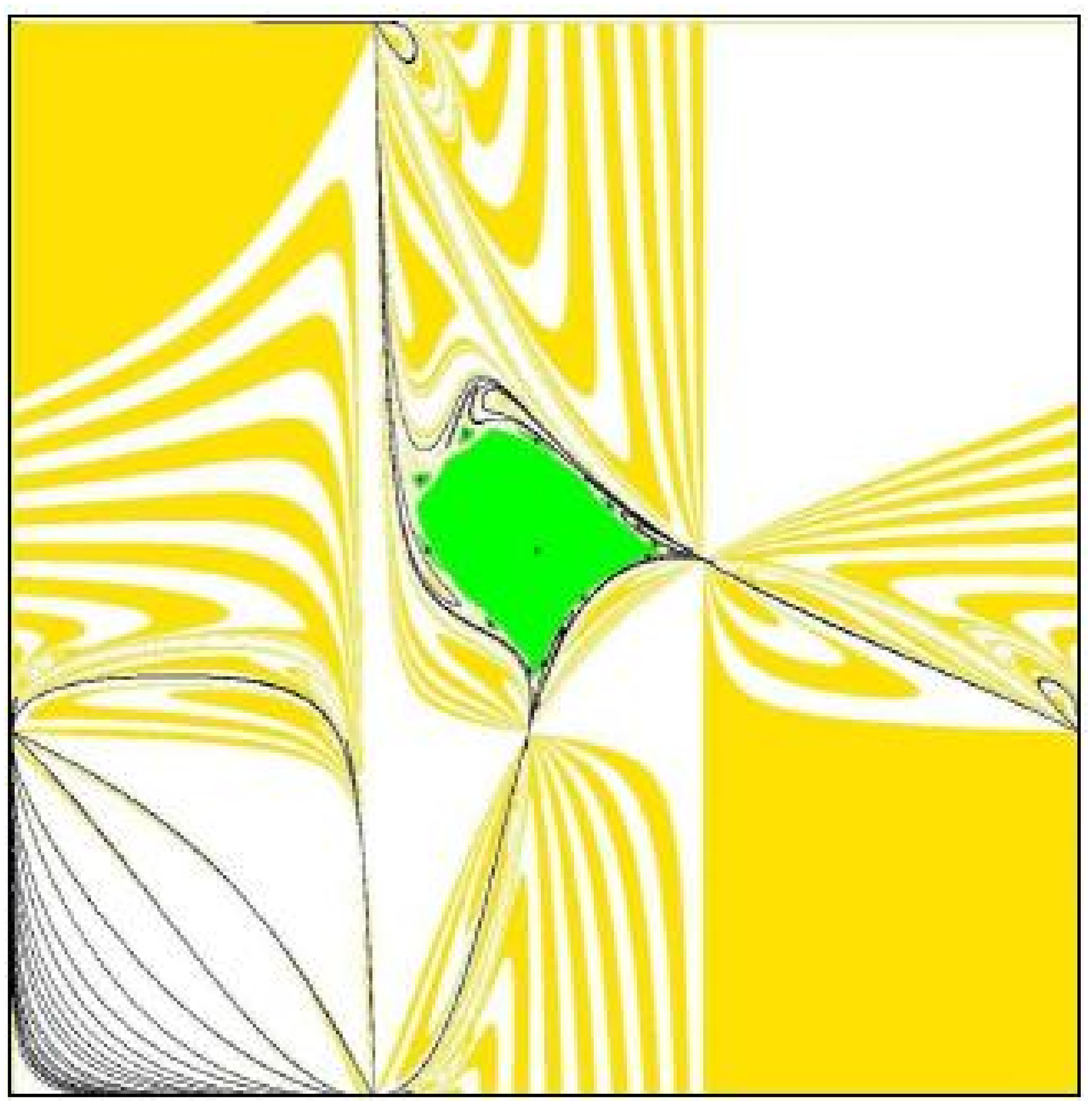}}
\resizebox{3in}{3in}{\includegraphics{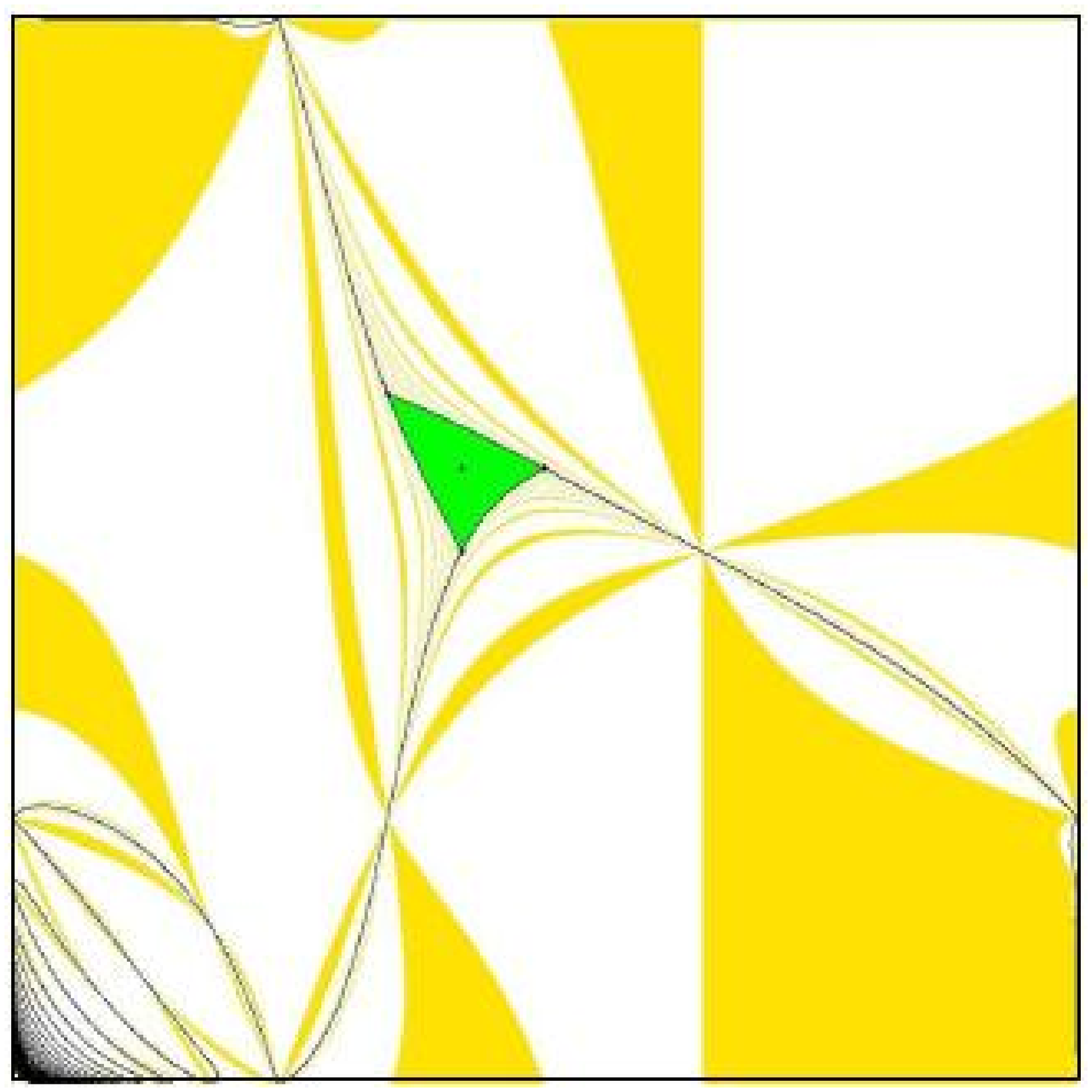}}
\end{center}
\caption{
\label{opimage}
Dynamics of $f$ for parameter values $a=1.1$ (left) and $a=2$ (right).  
Forward orbits of white points escape to infinity traveling up and to the 
right.  Yellow points escape by alternating between the bottom
right and upper left corners.  The black curves are stable manifolds of a 
saddle three cycle (the vertices of the `triangle').  The green regions 
consist of points whose orbits are bounded in forward and backward time.  
Reflection about the line $y=-x$ corresponds to replacing $f$ by $f^{-1}$.  It
leaves the green region invariant and exchanges stable and unstable manifolds.}
\end{figure} 

In \cite{BeDi05} we studied the dynamics of $f_a$ for the parameter region 
$a<0$, $a\ne-1$.  In this case, $\basin^+$ and $\basin^-$ are dense in $\R^2$;
and the complement $\R^2-\basin^+\cup\basin^-$ is a non-compact set on 
which the action of $f_a$ is very nearly hyperbolic and essentially 
conjugate to the golden mean subshift.  In particular, $f_a$ is topologically 
mixing on $\R^2-\basin^+\cup\basin^-$, and most points therein have
unbounded orbits.  The situation is quite different when $a>1$.

\begin{thm}
\label{mainthm1} 
If $a>1$, then $\R^2-\basin^+\cup \basin^- = K$.  Moreover, the set $K$ is 
compact in $\R^2$, contained in the square $[-a,1]\times [-1,a]$.
\end{thm}

The sets $\basin^+$ and $K$ are illustrated for typical parameter values
$a>1$ in Figure \ref{opimage}.   Note that all figures in this 
paper are drawn with $\R^2$ compactified as a torus $S^1\times S^1\cong
(\R\cup\{\infty\})\times (\R\cup\{\infty\})$.  The plane is parametrized so 
that infinity is visible, with top/bottom and left/right sides identified with 
the two circles at infinity.  All four corners of the square correspond to the
point $(\infty,\infty)$, which is a parabolic fixed point for $f_a$.

\begin{figure}
\begin{center}
\resizebox{3in}{3in}{\includegraphics{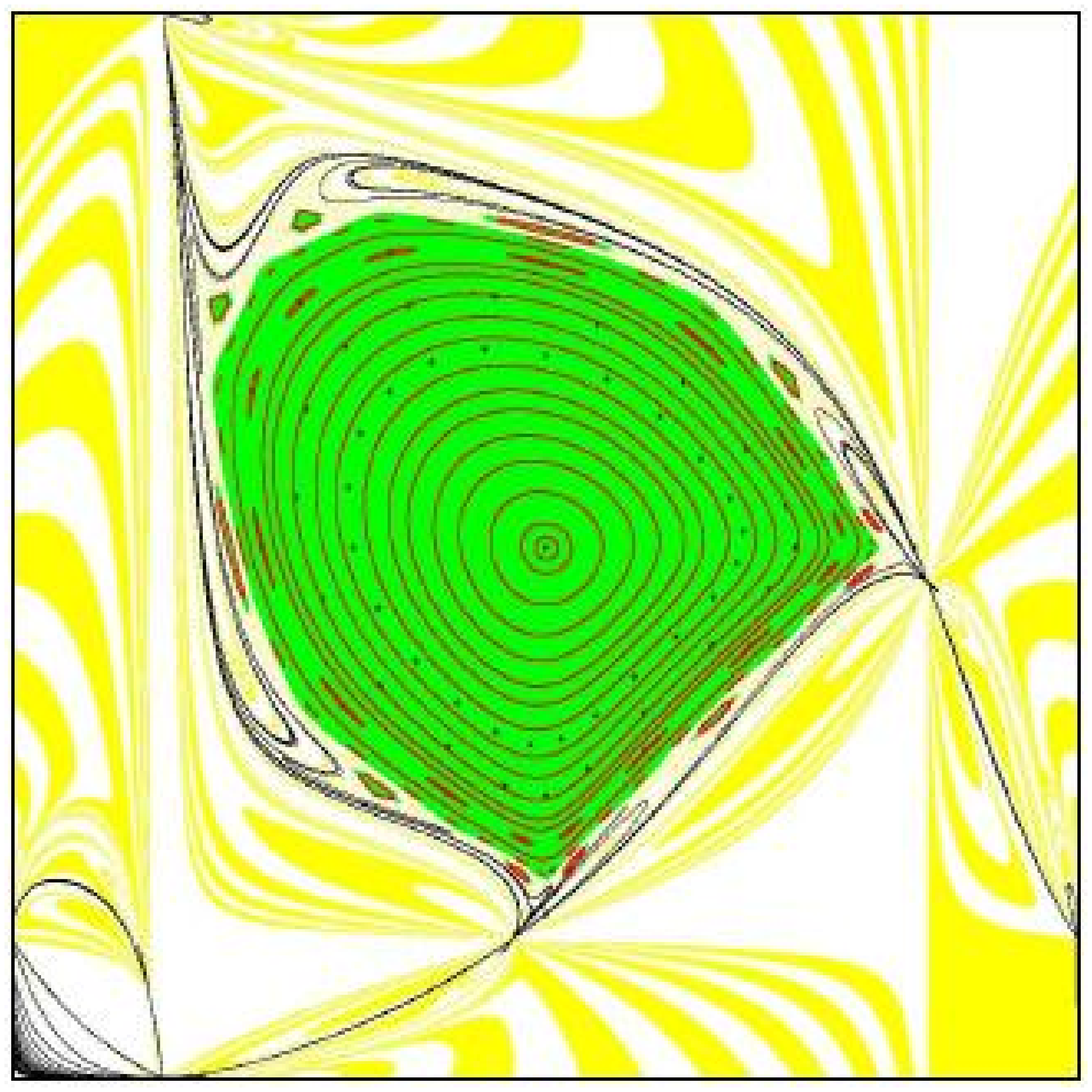}}
\end{center}
\caption{
\label{circles}
This is the left side of Figure \ref{opimage} redrawn with coordinates changed
to magnify the behavior near $p_{fix}$.  Shown are several 
invariant circles about $p_{fix}$, as well as some intervening
cycles and elliptic islands of large period.}
\end{figure} 

For each $a\neq -1$ there is a unique fixed point 
$p_{fix}=((1-a)/2,(a-1)/2)\in{\bf R}^2$.  For $a<0$, $p_{fix}$ is a saddle 
point, and for $a>0$, $p_{fix}$ is indifferent with $Df_a(p_{fix})$ conjugate
to a rotation.  For generic $a>0$ this rotation is irrational, and 
$f_a$ acts on a neighborhood of $p_{fix}$ as an area-preserving twist map with 
non-zero twist parameter (see Proposition \ref{bnfprop}).  As is shown in 
Figure 
\ref{circles}, $f_a$ exhibits KAM behavior for typical $a>0$.  In particular
$K$ has non-empty interior, and the restriction $f_a:K\self$ is neither 
topologically mixing nor hyperbolic.

\begin{figure}
\begin{center}
\resizebox{4in}{4in}{\includegraphics{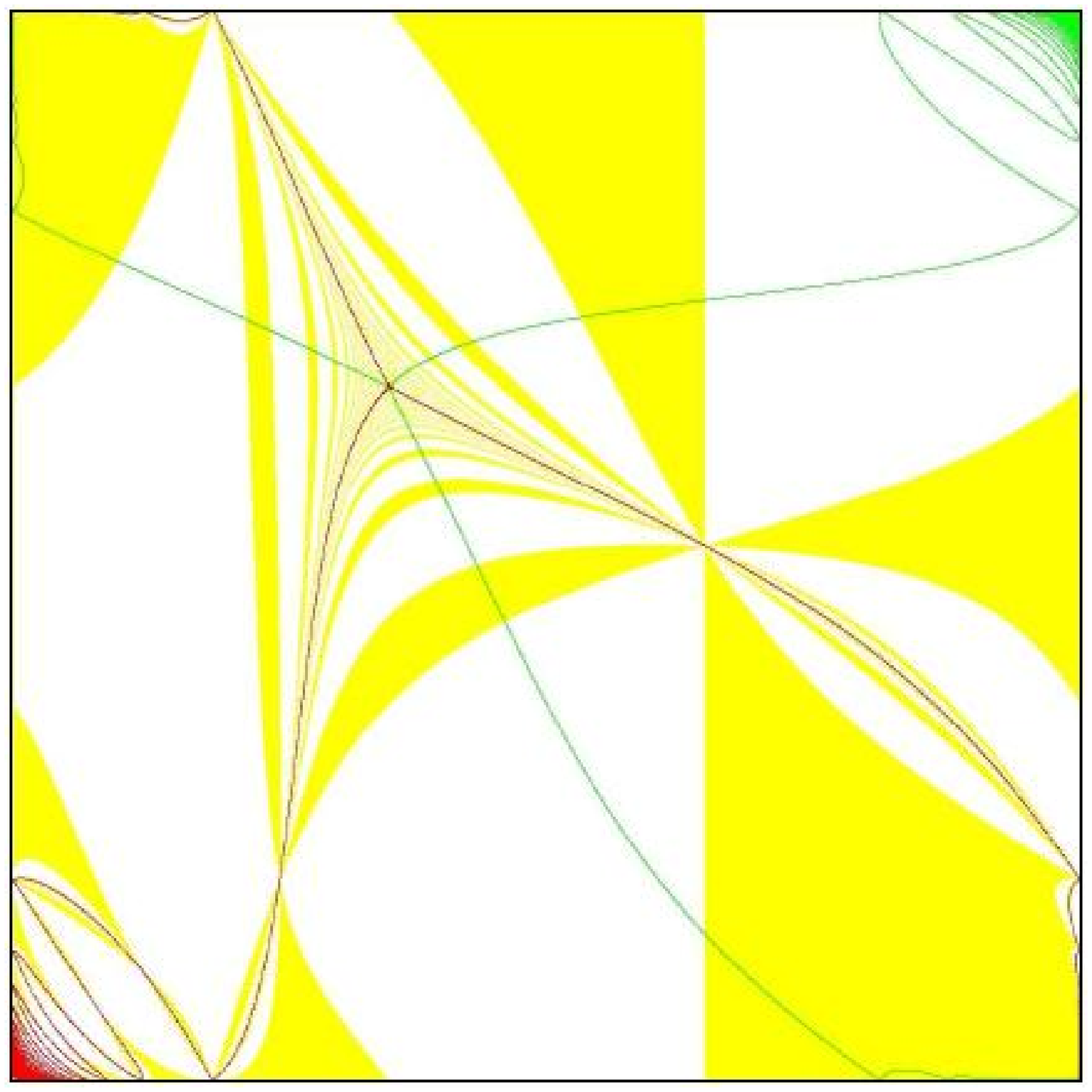}}
\end{center}
\caption{
\label{3}
Dynamics of $f$ when $a=3$.  The three cycle has now disappeared, and with it
the twistmap dynamics.  The indifferent fixed point is now parabolic, and its
stable and unstable sets are shown in red and green, respectively.}
\end{figure} 

It is known (see \cite{DiFa01} section 9) that $f_a$ is integrable for the 
parameter values $a=-1,0,\frac13,\frac12,1$.  In particular $f_a$ has
topological entropy zero for these parameters.   On the basis of computer 
experiments, \cite{AABM1,AABM3} conjectured that $a=3$ is the (unique) other 
parameter where 
entropy vanishes.   The map $f_3$ is not integrable, and in fact the 
complexified map $f_3:\C^2\self$ has entropy $\log{1+\sqrt5\over2}>0$ (see
\cite{BeDi05} and \cite{Du}).  
Nevertheless, we prove in Section \ref{aequals3} that not only does 
$f_3:\R^2\self$ have zero entropy; all points except $p_{fix}$ in $\R^2$ are 
transient: 

\begin{thm} 
\label{mainthm3} For $a=3$, $K=\{p_{fix}\}$.  The stable (and unstable) set of 
$p_{fix}$ consists of three analytic curves.  The three curves meet 
transversely at $p_{fix}$ and have no other pairwise intersections in 
${\bf R}^2$.  
\end{thm}

Thus, $p_{fix}$ is the only periodic and the only non-wandering point in
$\R^2$.  Figure \ref{3} illustrates this theorem.  As the figure makes
evident, the stable arcs for $p_{fix}$ intersect pairwise on a countable set 
at infinity.  This reflects the fact that $f$ is not a homeomorphism.

Often in this paper we will bound the number of intersections between two
real curves by computing the intersection number of their complexifications.
We will also use the special structure (see Figure \ref{partition}) of $f_a$
to help track the forward and backward images of curves.  The versatility
of these two techniques is seen from the fact that they apply in cases where
maps have maximal entropy (\cite{BeDi05,BeDi06}) and in the present situation
where we will show that $f_3$ has zero entropy. 


\section{Background}
\label{background}
\begin{figure}
\begin{center}
\resizebox{2.6in}{2.6in}{\includegraphics{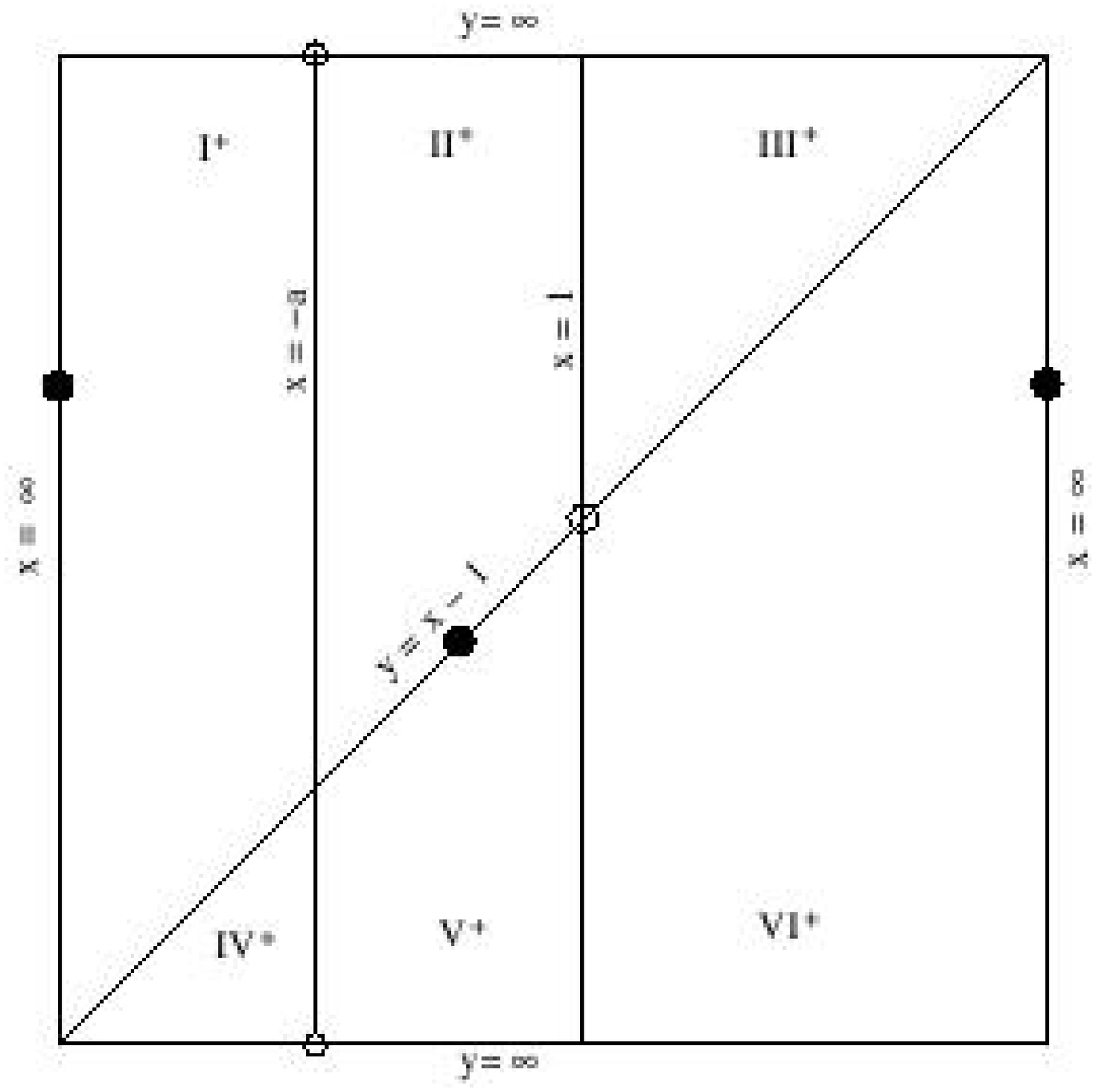}}
\hspace*{.3in}
\resizebox{2.6in}{2.6in}{\includegraphics{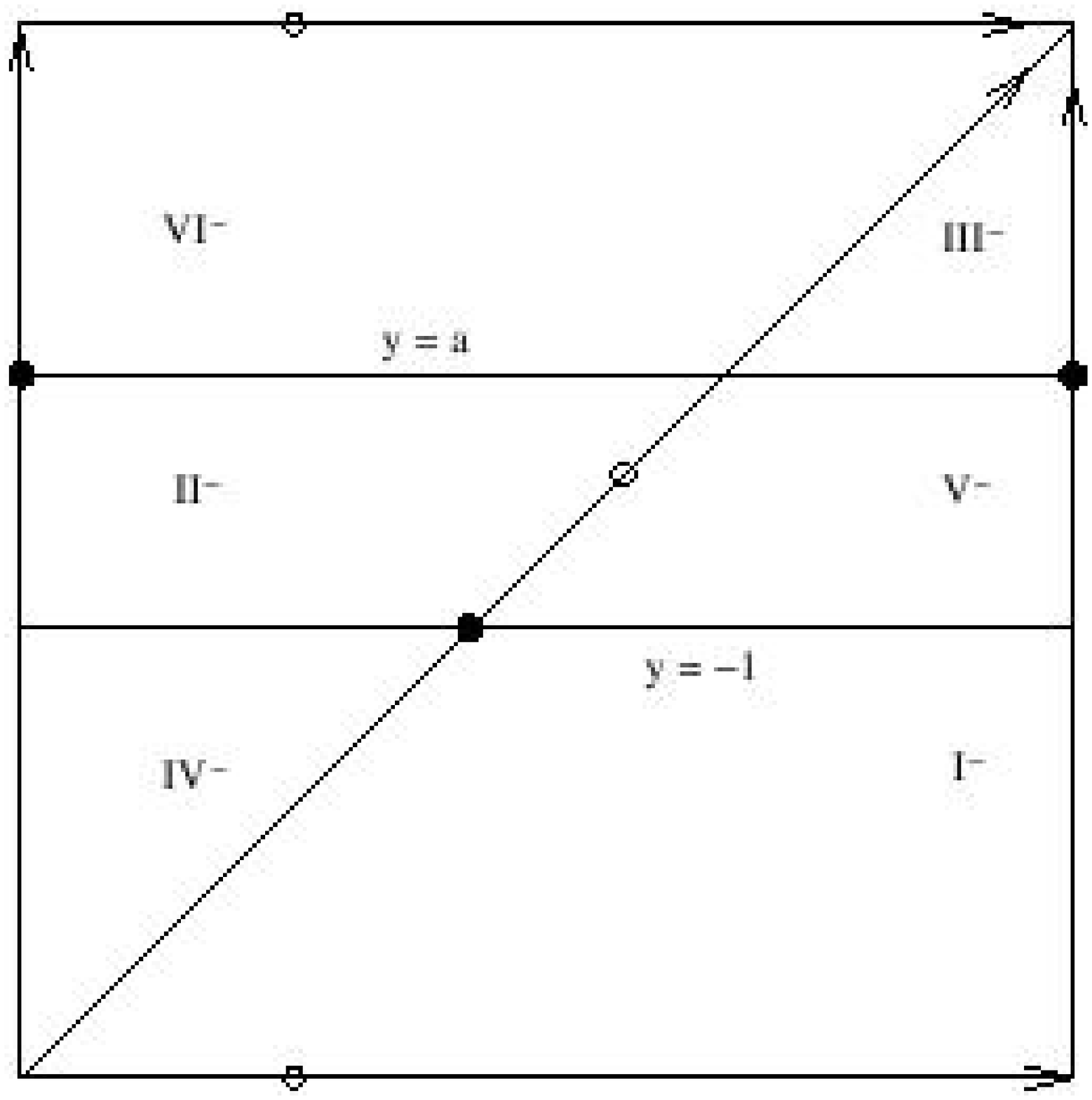}}
\end{center}
\caption{
\label{partition}
Partition of $\R^2$ by invariant curves and critical sets of $f$ (left side)
and $f^{-1}$ (right side).  Points of indeterminacy for $f$ appear as hollow
circles, whereas those for $f^{-1}$ are shown as solid circles.  Arrows on
the right side indicate the direction in which $f^2$ translates points along
$\supp(\eta)$.  }
\end{figure}

Here we recall some basic facts about the family of maps \eqref{opf}.  Most of
these are discussed at greater length in \cite{BeDi05}.
The maps \eqref{opf} preserve the singular two form 
$$
\eta \eqdef \frac{dx\wedge dy}{y-x+1}.
$$ 
Each map $f=f_a$ is also reversible, which is to say equal to a composition of 
two involutions \cite{BHM,BHM2}.  Specifically, $f=\tau\circ \sigma$
$$
\tau(x,y) \eqdef \left(x\frac{a-y}{1+y},a-1-y\right), \quad 
\sigma(x,y) \eqdef (-y,-x).
$$
In particular, $f^{-1} = \sigma\circ \tau$ is conjugate to $f$ by either 
involution, a property that will allow us to infer much about $f^{-1}$ directly
from facts about $f$.

Though our goal is to understand the dynamics of $f$ acting on $\R^2$, it will
be convenient to extend the domain of $f$ first by complexification to all of 
$\C^2$ and then by compactifying each coordinate separately to 
$\cp^1\times \cp^1$.  This gives us a convenient way of keeping track of the
complexity of algebraic curves in $\cp^1\times\cp^1$.  Any such curve
$V$ is given as the zero set of a rational function 
$R:\cp^1\times\cp^1\to \cp^1$.  The \emph{bidegree} $(j,k)\in\N^2$ of $R$, 
obtained by taking the degrees of $R$ as a function of the first and second 
variables in turn, encodes the second homology class of $V$.  In particular,
if $V$ and $V'$ are curves with bidegrees $(j,k)$ and $(j',k')$ having no 
irreducible components in common, then the number of (complex) intersections,
counted with multiplicity, between $V$ and $V'$, is 
\begin{equation}
\label{intersection}
V\cdot V' = jk' + j'k.
\end{equation}  
We will use
this fact at key points as a convenient upper bound for the number of
\emph{real} intersections between two algebraic curves.

Bidegrees transform linearly under our maps.  Taking
$f^* V$ to be the zero set of $R\circ f$, we have
\begin{equation}
\label{picactsbackward}
\operatorname{bideg} f^*V = 
\left(\begin{matrix} 1 & 1 \\ 1 & 0 \end{matrix}\right)
\operatorname{bideg} V.
\end{equation}
Similarly, $f_*V$ (the zero set of $R\circ f^{-1}$) has bidegree given by
\begin{equation}
\label{picactsforward}
\operatorname{bideg} f_*V = 
\left(\begin{matrix} 0 & 1 \\ 1 & 1 \end{matrix}\right)
\operatorname{bideg} V.
\end{equation}

The divisor $(\eta)$ of $\eta$, regarded as a meromorphic two form on
$\cp^1\times\cp^1$, is supported on the three lines $\{x=\infty\}$, 
$\{y=\infty\}$, $\{y=x-1\}$ where $\eta$ has simple poles.  It follows from 
$f^*\eta = \eta$ that $\supp(\eta)$ is invariant under $f$.
Specifically, $f$ interchanges the lines at infinity according to 
$$
(\infty,y) \mapsto (y,\infty) \mapsto (\infty,y+a-1)
$$ 
and maps $\{y=x-1\}$ to itself by $(x,x-1)\mapsto(x+a,x+a-1)$.  Thus
$f^2$ acts by translation on each of the three lines separately.  The
directions of the translations divide parameter space into three intervals:
$(-\infty,0)$, $(0,1)$, and $(1,\infty)$.  In particular, the directions are
the same for all $a\in (1,\infty)$, which is the range that concerns us here.

It should be stressed that $f:\cp^1\times \cp^1\self$ is \emph{not} a
diffeomorphism, nor indeed even continuous at all points.   In particular,
the critical set $C(f)$ consists of two lines 
$\{x=-a\}$ and $\{x=1\}$, which are  
mapped by $f$ to points $(0,1)$ and $(\infty,a)$, respectively.  
Applying the involution $\sigma$, one finds that $\{y=-1\}$ and $\{y=a\}$ are 
critical for $f^{-1}$,
mapping backward to $(-a,\infty)$ and $(-1,0)$, respectively.  Clearly,
$f$ cannot be defined continuously at the latter two points.  Hence we call
each a \emph{point of indeterminacy} and refer to 
$I(f) = \{(-a,\infty),(-1,0)\}$ as the \emph{indeterminacy set} of $f$.
Likewise, $I(f^{-1}) = \{(0,1),(\infty,a)\}$.  For convenience, we let
\begin{eqnarray*}
I^\infty(f) & \eqdef & \bigcup_{n\in\N} I(f^n) = \bigcup_{n\in\N} f^{-n} I(f)\\
C^\infty(f) & \eqdef & \bigcup_{n\in\N} C(f^n) = 
                       \bigcup_{n\in\N} f^{-n} C(f)
\end{eqnarray*}
denote the set of all points which are indeterminate/critical for high enough 
forward iterates of $f$.  We point out that both $I^\infty(f)$ and 
$I^{\infty}(f^{-1})$ are contained in $\supp(\eta)$.  Since $I(f)$ is
contained in $\supp(\eta)$, it follows that $I^\infty(f)$ is a discrete subset
of $\supp(\eta)$ that accumulates only at $(\infty,\infty)$.  We also observe 
that for the parameter
range $a>1$, we have $I^{\infty}(f)\cap I^\infty(f^{-1}) =
\emptyset$, a fact which will be useful to us below. 

The closure $\overline{\R^2}$ of the real points in $\cp^1\times\cp^1$ is just 
the torus $S^1\times S^1$.   The left side
of Figure \ref{partition} shows how the critical set of $f$ 
together with $\supp(\eta)$ partition $\overline{\R^2}$ into six open 
sets, which we have labeled $I^+$ to $VI^+$.  The right side shows
the partition by $\supp(\eta)$ and the critical set of $f^{-1}$.  As $f$
is birational, each piece of the left partition maps diffeomorphically onto
a piece of the right partition.  We have chosen labels on the right side so
that $I^+$ maps to $I^-$, etc.


\section{Trapping regions}
\label{trapping}
For the rest of this paper we assume $a>1$.
Figure \ref{partition} is useful for determining images and preimages of real 
curves by $f$.  In this section, we combine the information presented in the 
figure with intersection data gleaned from bidegrees to help identify two 
`trapping regions' through which orbits are forced to wander to infinity.  
Our first trapping region is the set
$$
T_0^+ \eqdef \{(x,y)\in\R^2: x>1, y>a\}
$$
of points lying to the right of the critical set of $f$ and above the critical 
set of $f^{-1}$.  

\begin{thm}
\label{trap0thm}
The region $T_0^+$ is forward invariant by $f$.  Forward orbits of points
in $T_0^+$ tend uniformly to $(\infty,\infty)$.
\end{thm}

\begin{proof}
Since $T_0^+\subset III^+\cup VI^+$ contains no critical or indeterminacy
points of $f$, we have that $f(T_0^+)$ is a connected open subset.  Since the 
left side of $T_0^+$ maps to a point, and 
$\partial T_0^+\cap I(f) = \emptyset$, we see that in fact $f(T_0^+)$ is
the region in $III^-\cup VI^-$ lying above and to the right of an 
arc $\gamma\subset f\{y=a\}$ that joins $(a,\infty) = f(\infty,a)$ to 
$(\infty,a) = f\{x=1\}$.  

By
\eqref{picactsforward}, the bidegree of $f\{y=a\} = f_*(y=a)$ is $(1,1)$.
Hence by \eqref{intersection}, we see that $f\{y=a\}$ intersects $y=a$ in
a single (a priori, possibly complex) point.  Since $(a,\infty)$ is one
such intersection, we conclude that $\gamma\cap\{y=a\}$ contains no points 
in $\R^2$.  That is, $\gamma$ lies above $\{y=a\}$.  Similarly, 
$\gamma\cap\{x=1\}$ contains at most one point.  However, because $a>1$, both 
endpoints of $\gamma$ have $x$ coordinates greater than $1$, and the number of
intersections between $\gamma$ and $\{x=1\}$ is therefore even.  We conclude 
that $\gamma\cap\{x=1\} = \emptyset$.  This proves $f(T_0^+)\subset T_0^+$.

\begin{figure}
\begin{center}
\resizebox{4in}{4in}{\includegraphics{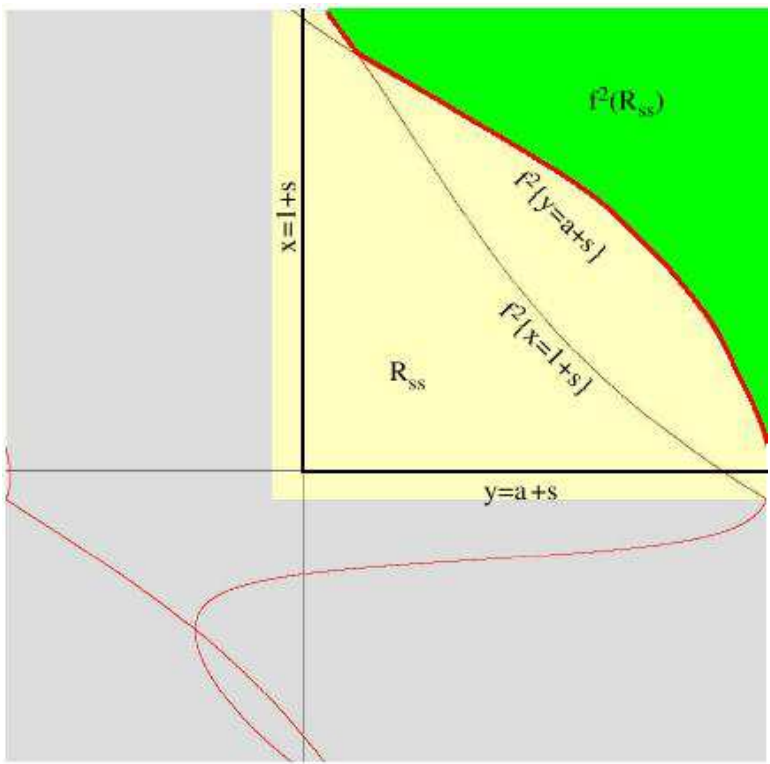}}
\end{center}
\caption{
\label{trap1}
First trapping region (shaded yellow/green) for $f$.  Also shown are the 
subregion
$R_{ss}$ (the square with lower left corner $(1+s,a+s)$ bounded by the thick 
black lines) and
$f^2(R_{ss})$ (green region) for $s=.3$.  The parameter value is $a=1.2$}.
\end{figure}

The same method further applies to show that points in $T_0^+$ have orbits 
tending uniformly to $(\infty,\infty)$.  If for $s,t>0$, we let 
$$
R_{st} = \{(x,y)\in \R^2:x>1+s,y>1+t\}\subset T_0^+,
$$
then arguments identical to those above combined with the fact that 
$f\{x=1+s\} = \{y=a+s\}$ suffice to establish
$$
f(R_{st}) \subset R_{t+a,s}.
$$
Hence $f^2(R_{ss}) \subset R_{s+a-1,s+a-1}$.  Since $a>1$ and the closures
$\overline{R_{ss}}\subset\overline{\R^2}$ decrease to the point
$(\infty,\infty)$ as $s\to\infty$, the proof is finished.
\end{proof}

Our second trapping region is more subtle.  In particular, it has two
connected components which are interchanged by $f$.  The first component is 
$$
A \eqdef
\{(x,y)\in\R^2:x>1,y<-x\}.  
$$
\begin{lem}
\label{trap2lem}
We have $A\cap f^{-1}(A) = \emptyset$ and $f(A)\subset f^{-1}(A)$.  
\end{lem}

\begin{figure}
\begin{center}
\resizebox{4in}{4in}{\includegraphics{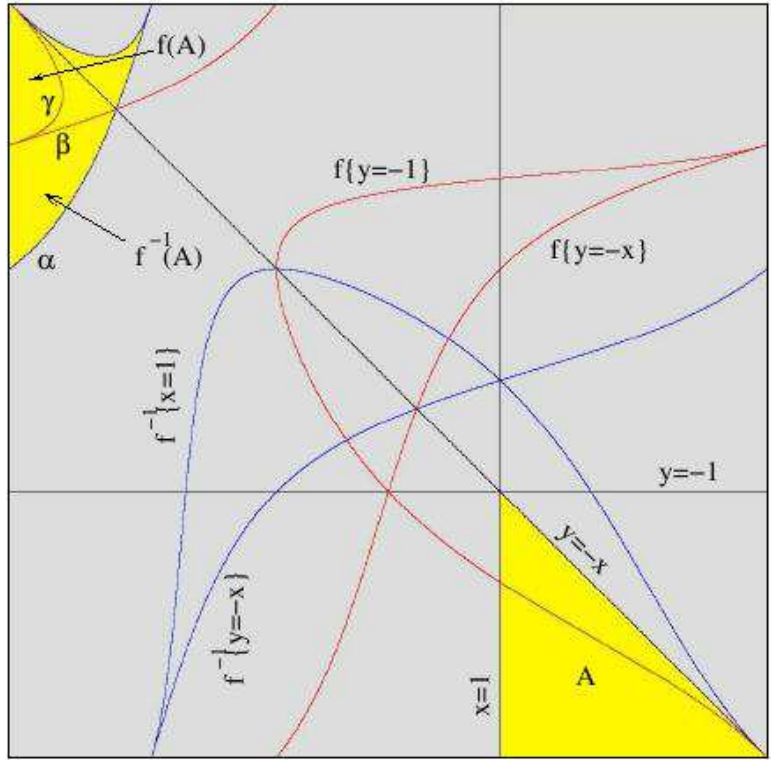}}
\end{center}
\caption{
\label{trap2}
The second trapping region for $(\infty,\infty)$, shown here in yellow
together with the curves used in the proof of Lemma \ref{trap2}.  Lines
used in the proof are black, images of lines are red, and preimages are
blue.  The parameter value is $a=3$.}
\end{figure}

\begin{proof}
The set $A$ lies entirely in region $I^-$ shown on the right side of Figure 
\ref{partition}, and therefore $f^{-1}(A)$ lies in region $I^+$ on the left
side.  In particular $A\cap f^{-1}(A) = \emptyset$.

To see that $f(A)\subset f^{-1}(A)$, we note that $A$ is bounded in 
$\overline{\R^2}$ by portions of three lines, one of which $\{x=1\}$ is
critical for $f$.  Moreover, $\overline{A}$ contains no point in $I(f)$
or $I(f^{-1})$ and $A$ itself avoids the critical sets of $f$ and $f^{-1}$.
Therefore, $f(A)$ is a connected open subset of $\R^2$, bounded on the left by 
the segment $\{(-\infty,y):y\geq a\}$ and the right by an arc $\gamma$ in 
$f\{y=-x\}$ that joins $(-\infty,a)$ to $(-\infty,\infty)$.  We claim that
$f(A)$ lies below $\{y=-x\}$ and above the arc $\beta$ in $f\{y=-1\}\cap VI^-$ 
that 
joins $(-\infty,a)$ to $(-a,\infty)$.  Both these claims can be verified in 
the same fashion; we give the details only for $\{y=-x\}$.

Clearly $(-\infty,a)\in \gamma$ lies below $\{y=-x\}$, and
$(-\infty,\infty)\in\gamma$ lies exactly on $\{y=-x\}$.  Hence to establish
that $f(A)$ is below $\{y=-x\}$ it suffices to show that $\gamma$ does not 
intersect $\{y=-x\}$ at any other point.  But the curve $f_*\{y=-x\}$
containing $\gamma$ has
bidegree $(1,2)$ and $\{y=-x\}$ has bidegree $(1,1)$, so there are a total
of three complex intersections between the two curves.  The intersection
at $(-\infty,\infty)$ is actually a tangency, which accounts for two of the
three intersections.  The fact that $\{y=-x\}$ joins the critical line
$\{x=-a\}$ to the critical line $\{x=1\}$ through regions $II^+$ and $V^+$
implies that $f\{y=-x\}$ joins through an arc passing right from $(1,0)$
to $(\infty,a)$.  This arc necessarily intersects $\{y=-x\}$ and accounts
for the third intersection.  Since there are no other intersections we
conclude that $\gamma$ is entirely below $\{y=-x\}$.

Having pinned down $f(A)$, we turn to $f^{-1}(A)$, which is bounded on the
left by $\{(-\infty,y):y\geq 1\}$, on the right by $f^{-1}\{x=1\}$ and from 
above by $f^{-1}\{y=-x\}$.  The top boundary component is in fact just
$\sigma(\gamma)$ and therefore lies above $\{y=-x\}$.  Therefore to show
$f(A)\subset f^{-1}(A)$, it suffices to show that $f(A)$ lies above the
portion $\alpha$ of $f^{-1}\{x=1\}$ bounding $f^{-1}(A)$.  For this, it 
suffices in turn to show simply that $\alpha$ lies below the arc $\beta$
described above.  This can also be accomplished with the method of the
previous paragraph and we spare the reader the details.
\end{proof}

We now define our second trapping region to be 
$$
T_1^+ \eqdef A \cup f^{-1}(A).  
$$
\begin{thm}
\label{trap1thm}
The region $T_1^+$ is forward invariant.  Forward orbits of points 
in $T_1^+$ tend uniformly to $(\infty,\infty)$, alternating between $A$
and $f^{-1}(A)$.
\end{thm}

\begin{proof}
Forward invariance of $T_1^+$ follows immediately from Lemma \ref{trap2lem}.
In order to show that points in $T_1^+$ have orbits tending uniformly to 
infinity, we consider diagonal lines
$$
L_t \eqdef \{y=t(x-1)\}, \quad L'_t \eqdef \{y=tx - 1\}
$$
passing through the points $(1,0) \in I(f)$ and $(0,-1)\in I(f^{-1})$, 
respectively.

\begin{lem}
The curve $f^2(L_t)\cap A$ (non-empty only for $t<-1$) lies strictly below 
$L'_t$.
\end{lem}

\begin{proof}
One computes easily that $f(L_t) = L'_{1/t}$ (we only include the
\emph{strict} transform of $L_t$ in the image $f(L_t)$; the line 
$\{y=a\} = f(1,0)$ appearing in $f_*L_t$ is omitted).  Because
$Df^2_{(\infty,\infty)} = \id$, it follows that $f^2(L_t) = f(L'_{1/t})$ is 
tangent to $L'_t$ at $(\infty,\infty)$.  Moreover, $f^2(L_t)\cap II^- = 
f(L'_{1/t}\cap II^+)$ joins $(0,-1)$ to $(a,0)$ and therefore intersects
$L'_t$ at $(0,-1)$.  This gives us a 
total of three points (counting multiplicity) in $f^2(L_t)\cap L_t' - A$.
On the other hand both $L'_t$ and $L_t$ have bidegree $(1,1)$ regardless of 
$t$, so from \eqref{intersection} and \eqref{picactsforward}, we find 
$L'_t\cdot f_* L'_{1/t}=3$.  It follows that 
$$
f^2(L_t)\cap L'_t\cap A = \emptyset.
$$  
Finally, since $f^2(L_t) \cap I^- = f(L'_{1/t}\cap I^+)$ joins $(0,-1)$ to 
$(\infty,\infty)$, we conclude that it lies strictly below $L'_t$.
\end{proof}

To complete the proof of Theorem \ref{trap1thm}, let $p\in T_1^+$ be given.
We may assume in fact that $p\in A$.  
Then each point $p_n \eqdef (x_n,y_n) \eqdef f^{2n}(p)$, $n\geq 0$ 
lies in $A$. Therefore $p_n \in L_{t_n}$ for some $t_n < -1$.  Because 
$L'_t$ is below $L_t$, the previous lemma implies that $(t_n)$ is a decreasing 
sequence. More precisely, $t_{n+1} < r t_n$, where $r=r(x_n)<1$ increases to
$1$ as $x_n\to\infty$.

Now if $(p_n)$ does not converge to $(\infty,\infty)$, we have a subsequence
$(p_{n_j})$ such that $x_{n_j} < M < \infty$.  The previous paragraph implies
that by further refining this subsequence, we may assume that 
$p_{n_j}\to (x,\infty)$ for some $x<\infty$.  This, however, contradicts the
facts that $f$ is continuous on $\overline{A}$ and that 
$f^n(x,\infty)\to (\infty,\infty)$.
\end{proof}

Let us define the forward basin $\basin^+$ of $(\infty,\infty)$ to be
the set of points $p\in\overline{\R^2}$ for which there exists a neighborhood 
$U\ni p$ such that $f^n|U$ is well-defined for all $n\in\N$ and converges 
uniformly to $(\infty,\infty)$ on $U$.  Our definition of $\basin^+$ here
differs slightly from the one given in the introduction in that we now allow
points at infinity.  Note that $(\infty,\infty)\notin
\basin^+$, because $I^\infty(f)$ accumulates at $(\infty,\infty)$.

\begin{thm}
\label{unbounded orbits}
The basin $\basin^+$ is a forward and backward invariant, connected and open 
set that contains all points in
$\supp(\eta)\cup C^\infty(f) - I^\infty(f) -\{(\infty,\infty)\}$.  Moreover, 
the following are equivalent
for a point $p\in\overline{\R^2}$.
\begin{itemize}
\item The forward orbit of $p$ is well-defined and unbounded.
\item $p\in \basin^+$.
\item $f^n(p)$ is in the interior of $\overline{T_0^+\cup T_1^+}$ for 
      $n\in\N$ large enough.
\end{itemize}
\end{thm}

We remark that in this context, it might be more relevant to consider 
orbits that accumulate on $\supp(\eta)$ rather than orbits which are
unbounded.  Theorem \ref{unbounded orbits} holds regardless.

\begin{proof}
The basin $\basin^+$ is open and invariant by definition.
Since $I^\infty(f)\subset \supp(\eta)$, we have that $\overline{I^\infty(f)}$
is discrete and accumulates only at $(\infty,\infty)$.  Therefore
connectedness of $\basin^+$ follows from invertibility of $f$.

Observe now that the $\overline{T_0^+ \cup T_1^+}$ contains a neighborhood
of any point $(t,\infty),(\infty,t), (t,t-1)\in \supp(\eta)$ for which
$t$ is large enough.  From this it follows easily that 
$\supp(\eta)-\overline{I^\infty(f)}\subset\basin^+$.  Since $C(f)-I(f)$ maps
to $I(f^{-1})$ and $I^\infty(f^{-1})\cap I^\infty(f) = \emptyset$, it further
follows that $C^\infty(f) - I^\infty(f)\subset\basin^+$.

The statements in the final assertion are listed from weakest to strongest, so
it suffices to prove the first implies the last.  So suppose that 
$p\in\overline{\R^2}$ is a point whose forward orbit is well-defined and 
unbounded.  If $(f^n(p))_{n\geq 0}$ accumulates at 
$q\in\supp(\eta)-I^\infty(f)$, then 
it also accumulates at every point in the forward orbit of $q$.  We observed
in the previous paragraph that $f^n(q)$ lies in the interior of
$\overline{T_0^+\cup T_1^+}$.  Therefore since $f^n$ is continuous on a
neighborhood of $q$, it follows that $f^m(p)\in \basin^+$ for some large $m$.
By invariance, we conclude $p\in\basin^+$, too.

The remaining possibility is that the forward orbit of $p$ accumulates at 
infinity only at points in $I^\infty(f)$.  This is impossible for the
following reason.  Since $I^\infty(f)$ is disjoint from 
$I^\infty(f^{-1})$, and since $f$ and $f^{-1}$ are conjugate via the 
automorphism $\sigma$, we have that the backward iterates $(f^{-n})$ converge 
uniformly to $(\infty,\infty)$ on a neighborhood $U\ni I^\infty(f)$.  Thus
points in $U$ cannot recur to points in $I^\infty(f)$, and in particular,
the forward orbit of $p$ cannot accumulate on $I^\infty(f)$.
\end{proof}

Because $f$ and $f^{-1}$ are conjugate via $(x,y)\mapsto (-y,-x)$, we
immediately obtain trapping regions 
$$
T_0^- = \sigma(T_0^+), \quad
T_1^- = \sigma(T_1^+)
$$
for $f^{-1}$, for which the exact analogues of Theorems \ref{trap0thm},
\ref{trap1thm}, and \ref{unbounded orbits} hold.  In particular, the backward
basin $\basin^- = \sigma(\basin^+)$ of $(\infty,\infty)$ includes all points
in $C^\infty(f)$ and all points in 
$\supp(\eta)-I^\infty(f^{-1}) - \{(\infty,\infty)\}$.  Using once again the 
fact that $I^\infty(f)$ is disjoint from $I^\infty(f^{-1})$, we have

\begin{cor}
\label{unbounded set}
All points except $(\infty,\infty)$ with unbounded forward or backward orbits 
are wandering.  Such points include all of $\supp(\eta)$, $I^\infty(f)$,
$I^\infty(f^{-1})$, $C^\infty(f)$ and $C^\infty(f^{-1})$
\end{cor} 


\section{Points with bounded orbits}
\label{bounded}
In this section we turn our attention to the set
$$
K = \{p\in\R^2:(f^n(p))_{n\in\Z} \text{ is bounded} \}
$$
of points with bounded forward and backward orbits.  We begin by emphasizing
an immediate implication of Theorem \ref{unbounded orbits}.

\begin{cor}
\label{bounded set}
$K$ is a compact, totally invariant subset of $\R^2-\supp(\eta)$ containing no 
critical or indeterminate point for any iterate of $f$.  Any non-wandering
point in $\overline{R^2}$ except $(\infty,\infty)$ belongs to $K$. 
\end{cor}

\begin{figure}
\begin{center}
\resizebox{4in}{4in}{\includegraphics{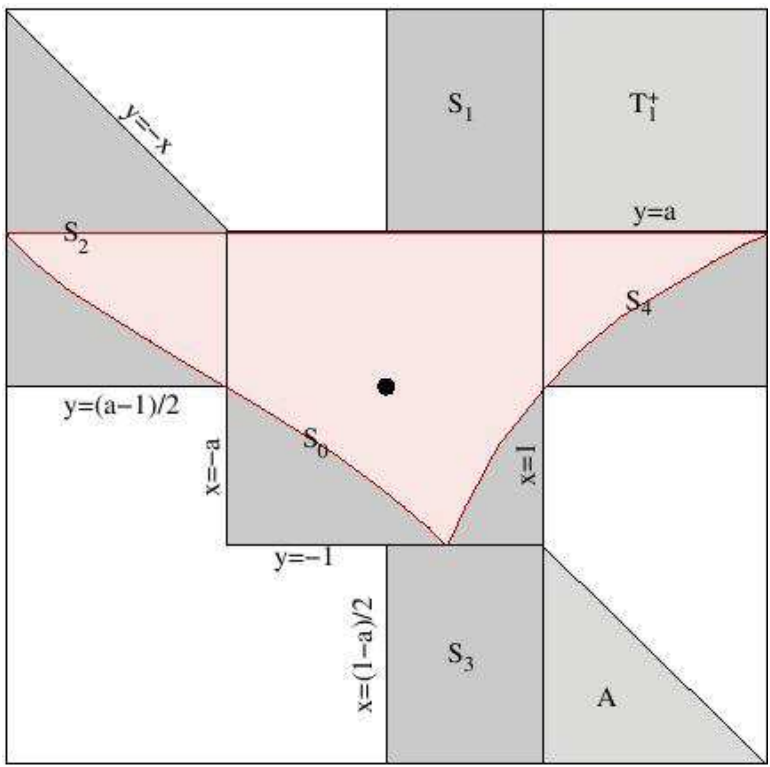}}
\end{center}
\caption{
\label{blades}
The regions $S_0$ through $S_4$ are quadrilaterals.  The dot in the center is 
the fixed
point, and the pink shaded region is $f(S_0)$.  The regions $T_0^+$
and $A$ are also shown.  Observe that if all shaded regions are reflected
about the line $y=-x$ (i.e. by the involution $\sigma$, then $S_0$ is sent
to itself and the images of the remaining regions 
exactly fill the complementary white regions.}
\end{figure}

To study $K$ we define several auxiliary subsets of $\R^2$.  Letting 
$(x_0,y_0) \eqdef (\frac{1-a}2,\frac{a-1}2)$ denote the coordinates of
the unique finite fixed point for $f$, we set 
$$
\begin{array}{ll}
S_0 \eqdef [-a,1]\times [-1,a],&
\quad S_1 \eqdef [x_0,1]\times [a,\infty], \\
S_2 \eqdef \{(x,y):x\leq -1, \frac12 (a-1) \leq y \leq -x\}, &
\quad S_3 \eqdef [x_0,1]\times [-\infty,-1],\\
S_4 \eqdef [1,\infty]\times [y_0,a].
\end{array}
$$
These sets are shown in Figure \ref{blades} with the image $f(S_0)$ 
superimposed.
Observe that $S_0$ is just the square in $\R^2$ cut out by the
critical sets of $f$ and $f^{-1}$, and that the other four regions, 
three rectangles and a trapezoid ($S_2$), are arranged around $S$ somewhat 
like `blades' on a fan.  Theorem \ref{bladethm} shows that 
$f$ `rotates' the blades counterclockwise.  It also shows that $S_0$ acts as
a kind of reverse trapping region for $K$ since a point in $S_0$ whose orbit
leaves $S_0$ cannot return again.

\begin{thm} 
\label{bladethm}
The sets $S_0,\dots, S_4$ satisfy the following:
\begin{enumerate}
\item $T_0^+, T_1^+, S_0,
\dots, S_4$, together with their images under the involution 
$\sigma$, cover $\R^2$;
\item $f(S_0)\subset S_0\cup S_2\cup S_4$;
\item $f(S_1)\subset S_2$;
\item $f(S_2)\subset S_3\cup A$;
\item $f(S_3)\subset S_4$;
\item $f(S_4)\subset S_1\cup T_0^+$;
\item any point $p\in S_1\cup \dots\cup S_4$ has a forward orbit tending to
$(\infty,\infty)$.
\end{enumerate}
\end{thm}
From this theorem, we have immediately that 
\begin{cor}
\label{trichotomy}
For every point $p\in\R^2$, we exactly one of the following is true.
\begin{itemize}
\item $f^n(p)\notin S_0$ for all $n\in\N$.
\item There exists $N\in\N$ such that $f^n(p)\notin S_0$ for $n<N$ and
$f^n(p)\in K$ for all $n\geq N$.
\item There exist $N\leq M\in \N$ such that $f^n(p)\in S_0$ for all
$N\leq n\leq M$ and $f^n(p)\notin S_0$ otherwise.
\end{itemize}
In the first two cases, $p$ has a forward orbit tending to $(\infty,\infty)$
through $T_0^+$ or $T_1^+$.  
\end{cor}

Since $f$ and $f^{-1}$ are conjugate by $\sigma$, we further have

\begin{cor}
\label{in the box}
$K\subset S_0$.
\end{cor}

Theorem \ref{unbounded orbits} and Corollary \ref{in the box} combine to yield 
Theorem \ref{mainthm1} in the introduction.
We spend the rest of this section proving Theorem \ref{bladethm}.  The first
assertion is obvious from Figure \ref{blades}.

\proofof{Assertion (ii)}
Since $S_0$ lies between the critical lines of $f$ in regions $II^+$ and $V^+$
(see Figure \ref{partition}), $f(S_0)$ is the closure of a connected open set
in regions $II^-$ and $V^-$.  Moreover, $\partial S_0$ is contains the point
$(1,0)\in I(f)$, which maps to $y=a$.  Hence $y=a$ constitutes the upper 
boundary of $S_0$.  

The vertical sides of $S_0$ are critical for $f$ and map to points, so $S_0$ 
is bounded below by the intersections $\alpha$ and $\beta$ of $f\{y=a\}$ and 
$f\{y=-1\}$, respectively, with $II^-\cup V^-$.  From Figure \ref{partition} 
one finds that $\alpha$ stretches right from $(-\infty,a)$ to $(0,-1)$ and
that $\beta$ stretches left from $(0,-1)$ to $(\infty,a)$.  Moreover, one
computes from \ref{picactsforward} and \ref{intersection} that each of these 
curves
intersects each horizontal line no more than once.  Hence $\alpha$ has
non-positive slope at all points and $\beta$ has non-negative slope.  Since
$(-a,y_0) = f(x_0,a)\in \alpha$ and 
$(1,y_0) = f(x_0,-1)\in\beta$, it follows that $f(S_0)-S_0$ is
contained in $S_2$ and $S_4$.
\qed

To verify the remaining assertions, it is useful to choose non-negative 
`adapted' coordinates on $S_1$ through $S_4$ as follows.  Each boundary 
$\partial S_j$, $1\leq j\leq 4$,
contains a segment $\ell$ of either the horizontal or vertical line passing 
through the fixed point $(x_0,y_0)$.  We choose coordinates 
$(x_j,y_j) = \psi_j(x,y)$ on $S_j$ so that 
$\ell\cap S_0$ becomes the origin (marked in Figure \ref{blades}, distance to 
$\ell$ becomes the $x_j$
coordinate and distance along $\ell$ to $\ell\cap S_0$ becomes the $y_j$ 
coordinate.  Thus, for example, we obtain coordinates
$
(x_1,y_1) = \psi_1(x,y) \eqdef (x-x_0,y-a)
$
on $S_1$ and coordinates
$
(x_2,y_2) = \psi_2(x,y) \eqdef (y-y_0,1-x)
$
on $S_2$, etc.

\proofof{assertions (iii) to (vi)}
To establish assertion (iii), we let 
$(x_1,y_1) = \psi_1(x,y)$ be the adapted 
coordinates of a point $(x,y)\in S_1$ and $(x_2,y_2) = \psi_2\circ f(x,y)$ the
coordinates of its image.  Then $0\leq x_1\leq \frac{a+1}{2}$ and $0\leq y_1$, 
and direct computation gives that 
$$
(x_2,y_2)\eqdef \psi_2\circ f \circ \psi_1^{-1}(x_1,y_1) = 
\left(x_1,\frac{4ax_1 + y_1+ay_1 +2x_1y_1}{1+a-2x_1}\right).
$$
Hence $(x_2,y_2)$ satisfy the same inequalities as $(x_1,y_1)$, and it follows
that $f(x,y)\in S_2$.  

Verifying assertion (iii) is more or less the same, though messier because
of $S_2$ is not a rectangle.  If $(x,y)\in S_2$, then 
$(x_2,y_2)\eqdef \psi(x,y)$ satisfies $0\leq x_2 \leq y_2+\frac{a+1}2$ and 
$0\leq y_2$.
One computes 
$$
(x_3,y_3) \eqdef \psi_3\circ f(x,y) 
= \left(\frac{-1 + a^2 + 2y_2(x_2 + a -1)}{2(y_2+a+1)},y_2\right).
$$ 
Since $a>1$, we see that $x_3,y_3\geq 0$.  So to complete the proof, it 
suffices to show that $x_3 \leq y_3 + \frac{a+1}{2}$:
$$
y_3 + \frac{a+1}{2} - x_3 = \frac{2+2a + y_2(5+a+2(y_2-x_2))}{2(y_2 + a +1)}
> \frac{y_2(5+a+(a+1))}{2(y_2+a+1)} > 0.
$$ 
We leave it to the reader to verify assertions (v) and (vi).
\qed

\begin{lem} 
\label{y gets bigger}
Let $p\in S_1\cup S_2\cup S_3\cup S_4$ be any point.  Then the
adapted $x$-coordinate of $f^2(p)$ is positive.  If, moreover, the adapted 
$y$-coordinate of $p$ is positive and $f^2(p) \in S_1\cup S_2\cup S_3\cup
S_4$, then the adapted $y$-coordinate
of $f(p)$ is larger than that of $p$. 
\end{lem}

\begin{proof}
In the same way we proved assertions (iii) through (vi) above, one may verify
that the restriction of $f$ to $S_1\cup S_3$ preserves the adapted
$x$-coordinate of a point, and the images of $S_2$ and $S_4$ do not contain
points with adapted $x$-coordinate equal to zero.  Hence $f^2(p)$ cannot
have adapted $x$-coordinate equal to zero.

Similarly, one finds that the restriction of $f$ to $S_2$ and $S_4$ preserves 
adapted $y$-coordinates, whereas $f$ increases the adapted $y$-coordinates
of points in $S_1$ and $S_3$.  For example, if $p = (x_1,y_1) \in S_1$, we
computed previously that $f(p)$ has adapted $y$-coordinate
$$
y_2 = \frac{4ax_1 + y_1+ay_1 +2x_1y_1}{1+a-2x_1} \geq y_1 + 
\frac{4ax_1}{1+a-2x_1} \geq y_1
$$
with equality throughout if and only if $x_1=0$.
\end{proof}

\proofof{assertion (vii)}
Suppose the assertion is false. Then by assertions (iii) through (vi) it 
follows that there is a point $p$ whose forward orbit is entirely contained
$S_1\cup S_2\cup S_3\cup S_4$.  Moreover, by Theorem \ref{unbounded orbits} the
forward orbit of $p$ must be bounded.  Therefore, we can choose 
$q\in S_1\cup S_2\cup S_3\cup S_4$ the be an accumulation point of
$(f^n(p))_{n\in\N}$ whose adapted $y$-coordinate is as large as possible.  

So on the one hand $f^2(q)$ is an accumulation point of the forward orbit of
$p$ and cannot, by definition of $q$, have adapted $y$-coordinate larger than
that of $q$.  But on the other hand, the first assertion of 
Lemma \ref{y gets bigger} implies that $q$ cannot have adapted $x$-coordinate
equal to zero, and therefore the second assertion tells us that
$f^2(q)$ must have larger adapted $y$-coordinate than $q$.  This contradiction
completes the proof.
\qed


\section{Behavior near the fixed point}
\label{twist}
When $a\geq 0$, the eigenvalues of $Df(p_{fix})$ are a complex 
conjugate pair $\lambda,\bar\lambda$ of modulus one, with
$$
\lambda = e^{i\gamma_0}=\frac{i+\sqrt a}{i-\sqrt a}.
$$  
Hence $Df(p_{fix})$ is conjugate to rotation by an angle $\gamma_0$ that
decreases from $0$ to $-\pi$ as $a$ increases from $0$ to $\infty$.
When $\gamma_0\notin \Q\pi$ is an irrational angle, it is classical that
$f$ can be put formally into Birkhoff normal form:

\begin{prop}
\label{bnfprop}
If $\lambda$ is not a root of unity, there is a formal change of coordinate
$z=x+iy$ in which $p_{fix} = 0$, and $f$ becomes 
$$
z \mapsto z\cdot e^{i(\gamma_0+\gamma_2|z|^2+\cdots)}
$$
where $\gamma_0$ is as above, and
$$
\gamma_2(a)=\frac{4(3a-1)}{\sqrt a (a-3)(1+a)^2}.
$$
\end{prop}

Hence the `twist' parameter $\gamma_2$ is well-defined and non-zero everywhere
except $a=\frac13$ and $a=3$, changing signs as it passes through these two 
parameters.  

For the proof of the proposition, we briefly explain how the Birkhoff normal 
form
and attendant change of coordinate are computed (see \cite{SM}, \S 23 for a
more complete explanation).
We start with a matrix $C$ whose columns are complex conjugates of each other 
and which satisfies 
$C^{-1}Df(p_{fix}) C=\mathrm{diag}(\lambda,\bar\lambda)$.  
Conjugating $C^{-1}\circ f_a\circ C$, we ``complexify'' the map $f_a$, using 
new variables $(s,t)=(z,\bar z)$, in which the map becomes
$$
f:  \ s\mapsto \lambda s + p(s,t),\ \ \  t\mapsto\bar\lambda t+q(s,t),
$$
and the power series coefficients of $q$ are the complex conjugates of the 
coefficients of $p$.  We look for a coordinate change 
$s=\phi(\xi,\eta)=\xi+\cdots$, $t=\psi(\xi,\eta)=\eta+\cdots$, such that the 
coefficients of $\phi$ and $\psi$ are complex conjugates, and which satisfies
$u\cdot \phi=p(\phi,\psi),$
for some $u(\xi,\eta)=\alpha_0+\alpha_2 \xi\eta+\alpha_4\xi^2\eta^2+\cdots$.  
Solving for the coefficients, we find 
$u={\rm exp}i(\gamma_0+\gamma_2 \xi\eta+\cdots)$, with $\gamma_0$ and 
$\gamma_2$ as above.   Returning to the variables $s=z$ and $t=\bar z$ gives 
the normal form.

For the rest of this section, we restrict our attention to the case $a=3$.  As 
was noted in \cite{BHM}, this is the parameter value for which a 3-cycle of 
saddle 
points coalesces with $p_{fix}$.  When $a=3$, we have $\gamma_0=-2\pi/3$.   
Translating coordinates so that $p_{fix}$ becomes the origin, we have
\begin{eqnarray*}
f^3(x,y) &=& (x,y) + Q+O(|(x,y)|^3) \\
Q =(Q_1,Q_2)&=& (x^2/2+xy+y^2,x^2+xy+y^2/2).
\end{eqnarray*}
Let us recall a general result of Hakim \cite{Hak} on the local structure of a 
holomorphic map which is tangent to the identity at a fixed point.  A vector 
$v$ is said to be characteristic if $Qv$ is a multiple of $v$.  The 
characteristic vectors $v$ for the map $f_a$ are $(1,-\frac12)$, $(1,1)$, 
and $(1,2)$.  Hakim \cite{Hak} shows 
that for each characteristic $v$ with $Qv\neq 0$, there is a holomorphic 
embedding 
$\varphi_v:\Delta\to{\bf C}^2$ which extends continuously to $\bar\Delta$ and 
such that $\varphi_v(1)=(0,0)$, and the disk $\varphi_v(\Delta)$ is tangent 
to $v$ at $(0,0)$.   Further, $f(\varphi_v(\Delta))\subset\varphi_v(\Delta)$, 
and for every $z\in\varphi_v(\Delta)$, $\lim_{n\to\infty}f^nz=(0,0)$.  
The disk $\varphi_v(\Delta)$ is a $\C^2$ analogue of an `attracting petal' at 
the origin. Applied 
to our real function $f_a$, this means that for each of the three vectors $v$, 
there is a real analytic ``stable arc'' $\gamma^s_v\subset{\bf R}^2$, ending 
at $(0,0)$ and tangent to $v$.  Considering $f^{-1}$, we have an ``unstable 
arc''  $\gamma^u_v$ approaching $(0,0)$ tangent to $-v$.  These two arcs fit 
together to make a $C^1$ curve, but they are not in general analytic 
continuations of each other.

Let us set $r(1,y)=Q_2(1,y)/Q_1(1,y)$ and write $v=(1,\eta)$.  Then 
$a(v):= r'(\eta)/Q_1(1,\eta)$ is an invariant of the map at the fixed point.  
If we make a linear change of coordinates so that the characteristic vector 
$v=(1,0)$ points in the direction of the $x$-axis, then we may rewrite $f$ in 
local coordinates as
$$
(x,u)\mapsto (x-x^2+O(|u|x^2, x^3),u(1-ax)+O(|u|^2 x, |u|x^2)),    \eqno(6)
$$
(see \cite{Hak}).  For the function (5), we find that $a(v)=-3$ for all three 
characteristic vectors.  We conclude that the stable arcs $\gamma^s_v$ are 
weakly repelling in the normal direction.  

Another approach, which was carried out in \cite{AABM3}, is to find the formal
power series expansion of a uniformization of $\gamma^s_v$ at $p_{fix}$.


\section{The case $a=3$}
\label{aequals3}
In this section we continue with the assumption $a=3$,
our aim being to give a global treatment of the stable arcs for $p_{fix}$.
The basic idea here is that the behavior of $f|_{S_0}$ is controlled by 
invariant cone fields on a small punctured neighborhood of $p_{fix}$. 
 
In order to proceed, we fix some notation.  Let 
$\phi(x,y) = \phi_0(x,y) = x+1$ and for each $j\in\Z$, let 
$\phi_j(x,y)=\phi\circ f^{-j}(x,y)$.  Then $\phi_{j+k} = \phi_j\circ f^{-k}$
and in particular, $\phi_j(p_{fix}) = \phi_0\circ f^{-j}(p_{fix}) = \phi_0(p_{fix}) = 0$ for every $j\in\Z$.  We will be particularly concerned with the 
cases $j=-1,0,1,2$ and observe for now that the level set $\{\phi_j=s\}$ is 
\begin{itemize}
\item a horizontal line $\{x=-1+s\}$ when $j=0$;
\item a vertical line $\{y=1+s\}$ when $j=1$;
\item a hyperbola with asymptotes $\{x=-3\}$, $\{y=-1+s\}$ when $j=-1$; and
\item a hyperbola with asymptotes $\{x=1-s\}$, $\{y=3\}$ when $j=2$.
\end{itemize}
Most of our analysis will turn on the interaction between level
sets of $\phi_{-1}$ and $\phi_3$. 

\begin{prop}
\label{levelcurves}
$\{\phi_2=0\}$ is tangent to $\{\phi_{-1}=0\}$ at $p_{fix}$.  Moreover
$S_0\cap \{\phi_2>0\}\subset S_0\cap\{\phi_{-1}>0\}$.
\end{prop}

\begin{proof}
The first assertion is a consequence of the facts that 
$f^3\{\phi_{-1}=0\} = \{\phi_2=0\}$ and that $Df^3$ is the identity at 
$p_{fix}$.  Since both curves in question are hyperbolas with horizontal
and vertical asymptotes, it follows that $p_{fix}$ is the \emph{only} point
where the curves meet.  

The asymptotes of $\{\phi_{-1}=0\}$ are the bottom and left sides of $S_0$,
whereas those of $\{\phi_2=0\}$ are the top and right sides.  Therefore
the zero level set of $\phi_{-1}$ intersects $S_0$ in a connected, concave up 
arc; and the zero level set of $\phi_2$ meets $S_0$ in a connected concave down
arc.  It follows that the first level set is above and to the right of the
second.  Finally, direct computation shows that $\phi_{-1}$ and $\phi_2$
are positive at the lower left corner of $S_0$.  This proves the second assertion
in the proposition.
\end{proof}

Using the functions $\phi_j$, $j=-1,0,1$, we define `unstable wedges' 
$W_j\subset S_0$ emanating from $p_{fix}$:
\begin{eqnarray*}
W_0 &=& \{(x,y)\in S_0: 0 \leq \phi_{-1},\phi_0\} \\
W_1 &=& \{(x,y)\in S_0: 0 \leq \phi_0,\phi_1 \}\\
W_2 &=& \{(x,y)\in S_0: 0 \leq \phi_{-1},\phi_1\}
\end{eqnarray*}

\begin{prop}
\label{wedgescontract}
We have $f(W_0)\cap S_0 \subset W_1$, $f(W_1)\cap K \subset W_2$, 
and $f(W_2\cap S_0)\subset W_0$.  
\end{prop}

\begin{figure}
\begin{center}
\resizebox{4in}{4in}{\includegraphics{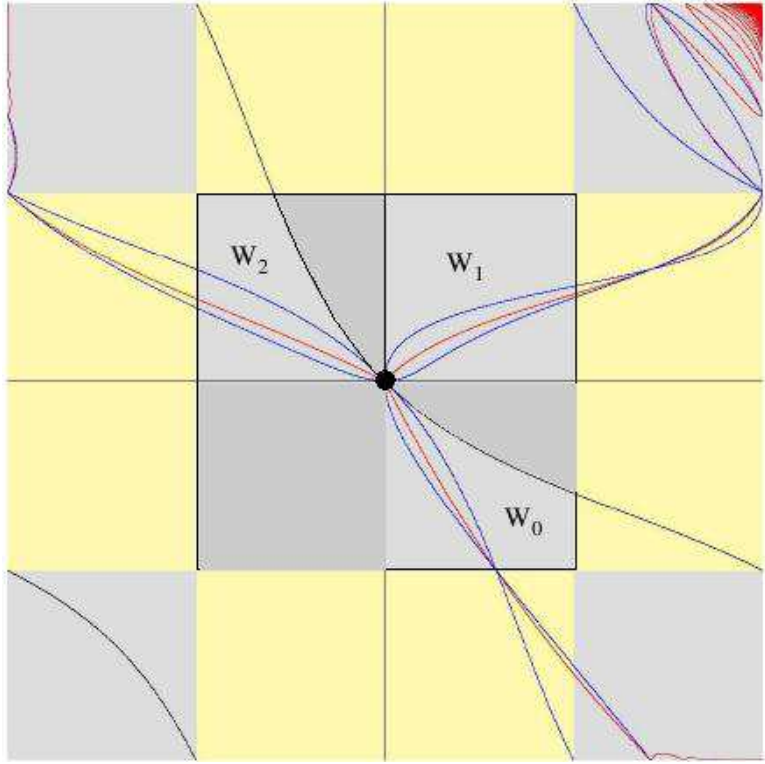}}
\end{center}
\caption{
\label{wedges}
Unstable wedges.  The square in the center is $S_0$, with the fixed point in
the middle.  Boundaries of the wedges are the black curves.  The images of
the wedges under $f^3$ are bounded by the blue curves.  The unstable curves
for the fixed point are red.  
}
\end{figure} 

\begin{proof}
It is immediate from definitions that 
$f(W_0)\cap S_0 \subset \{\phi_0,\phi_1\geq 0\}\cap S_0 = W_1$.

Likewise, 
$$
f(W_1)\cap S_0 \subset \{\phi_1,\phi_2\geq 0\}\cap S_0 \subset
\{\phi_1,\phi_{-1}\geq 0\}\cap S_0 = W_2.
$$
The second inclusion follows from the Proposition \ref{levelcurves}.
Similar reasoning shows that $f(W_2)\cap S_0 \subset W_0$.
\end{proof}

\begin{lem}
\label{comparison}
For each neighborhood $U$ of $p_{fix}$, there exists $m=m(U)>0$ such that
$$
\phi_{-1}(p)+m \leq \phi_2(p) \leq -m.
$$
for every $p\in W_0-U$.
\end{lem}

\begin{proof}
Given $p\in W_1$, we write $p=(-1+s,1+t)$ where $0\leq s,t \leq 2$ and
compute
$$
\phi_2(p) = \frac{-2s-2t-st}{t+2} = \frac{s(-2+t/2)+t(-2+s/2)}{t+2}
          < -\frac{s+t}{t+2} \leq -\frac{s+t}{4};
$$
and
$$
\phi_2(p) - \phi_{-1}(p) = \frac{2(s^2+t^2) + st(8+s-t)}{(2-s)(t+2)}
                      \geq \frac{2(s^2+t^2)}{4}
$$
Since the quantities $s^2+t^2$ and $s+t$ are both bounded below by positive 
constants on $W_0-U$, the lemma follows. 
\end{proof}

\begin{prop}
\label{allgo}
Let $U$ be any neighborhood of $p_{fix}$.  Then there exists $N>0$
such that $f^n(p)\notin S_0$ for every $p\in W_0\cup W_1\cup W_2-U$ and every 
$n\geq N$.  
\end{prop}

\begin{proof}
By Proposition \ref{wedgescontract}, we can assume that the point
$p$ in the statement of the lemma lies in $W_0-U$.
Since $\phi_2$ is continuous on $W_0$ with $\phi_2(p_{fix})=0$ and 
$\phi_2<0$ elsewhere on $W_0$, we may assume that $U$ is of the form
$\{\phi_2>-\epsilon\}$ for some $\epsilon>0$.  

Suppose for the moment 
that $f^3(p)\in S_0$.  Then Proposition \ref{wedgescontract} implies
that $f^3(p)\in W_0$.  By Lemma \ref{comparison}, we then have
$$
\phi_2\circ f^3(p) - \phi_2(p) = \phi_{-1}(p) - \phi_2(p) < -m(U).
$$
In particular, $f^3(p)\in W_0-U$.  Repeating this reasoning, we find
that if $p\in W_0-U$ and $f^{3j}(p)\in S_0$ for $j=1,\dots,J$, then 
$$
\phi_2\circ f^{3J}(p) < -Jm(U).
$$
On the other hand $\phi_2$ is bounded below on $W_0$ (e.g. by $-4$).
So if $J>4/m(U)$, it follows that $f^{3j}(p)\notin S_0$ for some $j\leq J$.
By Corollary \ref{trichotomy}, we conclude that $f^n(p)\notin S_0$ for all
$n\geq 3j$.
\end{proof}

Let us define the local stable set
of $p_{fix}$ to be
$$
W^s_{loc}(p_{fix}) \eqdef \{p\in\R^2:f^n(p)\in S_0 \text{ for all } n\in\N \text{
  and } f^n(p)\to p_{fix}\},
$$
with the local unstable set $W^u_{loc}(p_{fix})$ defined analogously.

\begin{cor}
\label{local stable}
The wedges $W_0,W_1,W_2$ meet $W^s_{loc}(p_{fix})$ only at $p_{fix}$, whereas 
they entirely contain $W^u_{loc}(p_{fix})$.  Thus
\begin{equation}
\label{equivalents}
W^s_{loc}(p_{fix}) = W^s(p_{fix})\cap S_0 
= \{p\in \R^2:f^n(p)\in S_0\text{ for all } n\in\N\},
\end{equation}
and similarly for $W^u_{loc}$.
\end{cor}

\begin{proof}
The first assertion about $W^s_{loc}(p_{fix})$ follows from Propositions 
\ref{wedgescontract} and \ref{allgo}.  Reversibility of $f$ then implies the 
first assertion about $W^u_{loc}(p_{fix})$ is disjoint from the sets
$\sigma(W_0),\sigma(W_1),\sigma(W_2)$.

We will complete the proof of the first assertion about $W_0$ 
$$
S_0 \subset \bigcup_{j=0}^2 W_j \cup \sigma(W_j).
$$
Observe that $\phi_0\circ\sigma = -\phi_1$.  Hence, 
$$
\phi_{-1}\circ\sigma = \phi_0\circ f \circ\sigma = -\phi_1\circ f^{-1} = 
-\phi_2.
$$  Thus $\sigma(W_0) = \{\phi_1,\phi_2 \leq 0\}$, etc.  The
desired inclusion is therefore an immediate consequence of Proposition 
\ref{levelcurves}.

The first equality in \eqref{equivalents} follows from Corollary 
\ref{trichotomy},
and the second from the analogue of Proposition \ref{allgo} for $f^{-1}$. 
\end{proof}

Each of the wedges $W_j$ admits an `obvious' coordinate system
identifying a neighborhood of $p_{fix}$ in $W_j$ with a neighborhood of 
$(0,0)$ in $\R^2_+ \eqdef \{x,y\geq 0\}$.  Namely, we let 
$\Psi_0:W_0\to \R^2$ be given by $\Psi_0 = (\phi_{-1},\phi_0)$,
$\Psi_1:W_1\to \R^2$ by $\Psi_1 = (\phi_0,\phi_1)$, and 
$\Psi_1:W_1\to \R^2$ by $\Psi_2 = (\phi_1,\phi_{-1})$.  If $U\subset\R^2_+$
is a sufficiently small neighborhood of $(0,0)$, then in light of 
Proposition \ref{wedgescontract}, the maps 
$f_{01},f_{12},f_{20}:U\to\R^2_+$ given by 
$f_{ij} = \Psi_j\circ f\circ \Psi_i^{-1}$ are all well-defined.  In fact,
$f_{01}$ is just the identity map.  The other two are more complicated, but
in any case a straightforward computation shows

\begin{prop}
\label{Df}
If $x,y\geq 0$ are sufficiently small, the entries of the matrix
$Df_{ij}(x,y)$ are non-negative.
\end{prop}

The proposition tells us that there are forward invariant cone fields
defined near $p_{fix}$ on each of the $W_j$; i.e. the cones consisting of 
vectors with non-negative entries in the coordinate system $\Psi_j$.  In this
spirit, we will call a connected arc $\gamma:[0,1]\in W_j\cap \overline{U}$ 
\emph{admissible} if $\gamma(0)=p_{fix}$, $\gamma(1)\in\partial U$ and  
both coordinates of the function $\Psi_j\circ\gamma:[0,1]\to\R^2$ are 
non-decreasing.  We observe in particular, that each of the two pieces of 
$\partial W_j\cap \overline{U}$ are admissible arcs.

By Proposition \ref{allgo}, we can choose the neighborhood $U$ of $p_{fix}$
above so that its intersections with the wedges $W_j$ are `pushed out' by $f$.
That is, $f(W_j\cap bU)\cap U = \emptyset$.  This assumption and
Propositions \ref{wedgescontract} and \ref{Df} imply immediately that 

\begin{cor}
\label{cones are invariant} If $\gamma:[0,1]\to W_i$ is an admissible arc, then
(after reparametrizing) so is $f_{ij}\circ\gamma$.
\end{cor}

\begin{thm}
\label{local local stable}
For each $j=0,1,2$, the set $\bigcap_{n\in\N} f^{3n}(W_j)\cap U$ is an 
admissible arc.
\end{thm}

\begin{proof}
For any closed $E\subset \overline{U}$, we let $\Area(E)$ denote the area
with respect to the $f$-invariant two form $\eta$.  Since $U$ avoids the
poles of $\eta$, we have $\Area(U) <\infty$.  

By Proposition \ref{wedgescontract} and our choice of $U$,
the sets $f^{3n}(W_j)\cap \overline{U}$ decrease as $n$ increases.  In 
particular, Proposition \ref{allgo} tells us that 
$\Area(f^{3n}(W_j)\cap \overline{U})\to 0$.  Finally, for every $n\in\N$,
$U\cap \partial f^{3n}(W_j) = f^{3n}(\partial W_j) \cap U$ consists of two
admissible curves meeting at $p$.  

Admissibility implies that in the coordinates $\Psi_j$, the two curves 
bounding $f^{3n}(W_j)$ in $U$ are both graphs over the line $y=x$ of functions
with Lipschitz constant no larger than $1$.  The fact that the region between 
these two curves decreases as $n$ increases translates into the statement that 
the graphing functions are monotone in $n$.  The uniformly bounded Lipschitz 
constant together with the fact that the area between the curves is tending 
to zero, implies further that the graphing functions converge uniformly to
the same limiting function and that this limiting function also has 
Lipschitz constant no larger than $1$.
\end{proof}

Observe that by Corollary \ref{local stable}, 
we have just characterized the intersection of the 
unstable set of $p_{fix}$ with $U$.  Since $f$ and $f^{-1}$ are conjugate
via the involution $\sigma$, the theorem also characterizes the stable set
of $p_{fix}$.  Our final result, stated for the stable rather than the 
unstable set of $p_{fix}$, summarizes 
what we have established.  It is an immediate consequence of Proposition 
\ref{allgo}, and Corollary \ref{trichotomy} and Theorem \ref{local local 
stable}.

\begin{cor} 
\label{last}
When $a=3$ the set $W^s(p_{fix})\cap S_0$ consists of three 
arcs meeting transversely at $p_{fix}$, and it coincides with the set of 
points whose forward orbits are entirely contained in $S_0$.  All points not 
in $W^s(p_{fix})$ have forward orbits tending to infinity.  In particular, 
$K=\{p_{fix}\}$ 
\end{cor}

Theorem \ref{mainthm3} is an immediate consequence.


\nocite{GrHa}
\nocite{KaHa}
\nocite{AABM1}
\nocite{AABM2}
\nocite{AABM3}
\nocite{AABM4}
\nocite{AABM5}
\nocite{AABM6}
\nocite{AABM7}
\nocite{Gue}

\bibliographystyle{/home/jeff/tex/mjo}
\bibliography{refs}
\end{document}